\documentclass[11pt]{amsart}                

\usepackage[all]{xy}
\usepackage{graphicx}
\usepackage{epstopdf}
\usepackage{pb-diagram}
\usepackage{amssymb, amscd}
\usepackage{amsmath}
\usepackage{amsfonts}
\usepackage{pstricks}
\usepackage{pst-grad}
\usepackage{pst-plot}
\usepackage[tiling]{pst-fill}
\usepackage{graphics}
 \xyoption{arc}

\theoremstyle{plain}
\newtheorem{theorem}{Theorem}[section]
\newtheorem{lemma}[theorem]{Lemma}

\newtheorem*{theorem*}{}
\newtheorem{prop}[theorem]{Proposition}

\theoremstyle{definition}

\newtheorem{definition}[theorem]{Definition}
\newtheorem{example}[theorem]{Example}
\newtheorem{remark}[theorem]{Remark}

\newcommand{\F}{{\mathcal F}}
\newcommand{\G}{\mathcal G}
\newcommand{\cL}{\mathcal L}
\newcommand{\Hol}{\mbox{\rm Hol}}

\newcommand{\id}{{\scriptstyle id}}
\newcommand{\Pair}{\mbox{\rm Pair}}

\newcommand{\Z}{{\mathbb Z}}

\newcommand{\cJ}{{\mathcal J}}
\newcommand{\cK}{{\mathcal K}}
\newcommand{\I}{{\mathcal I}}
\newcommand{\cP}{{\mathcal P}}

\newcommand{\cT}{{\mathcal T}}

\newcommand{\X}{{\mathcal X}}
\newcommand{\Y}{{\mathcal Y}}

\newcommand{\BB}{\mathsf{B}}
\newcommand{\GG}{{\mathsf G}}

\newcommand{\DD}{{\mathsf D}}
\newcommand{\MM}{\mathsf{Gpd}}
\newcommand{\HH}{\mathsf{H}}

\newcommand{\simM}{~ {\sim}_{M} ~}
\newcommand{\simT}{~ {\sim}_{T} ~}
\newcommand{\al}{\alpha}
\newcommand{\e}{\epsilon}

\newcommand{\ti}{\tilde}
\newcommand{\w}{\omega}
\newcommand{\ri}{\Rightarrow}
\newcommand{\inv}{^{-1}}
\newcommand{\globe}{
\[
 \xy
  (-15,0)*{\bullet}="1";
  (0,0)*{\bullet}="2";
 "1";"2" **\crv{(-12,9) & (-3,9)};
  "1";"2" **\crv{(-12,-9) & (-3,-9)};
    (-7.5,6.75)*{\scriptstyle >}+(0,3)*{\scriptstyle \phi};
    (-7.5,-6.75)*{\scriptstyle >}+(0,-3)*{\scriptstyle \psi};
    (-17.5,0)*{\scriptstyle \cK};
  (2.5,0)*{\scriptstyle \G};
  {\ar@{=>}_{\scriptstyle a}(-7.5,3)*{};(-7.5,-3)*{}} ;
 \endxy 
 \]
}
\newcommand{\vglobe}{
 \[
 \xy
  (-8,0)*{\bullet}="1";
  (8,0)*{\bullet}="2";
 "1";"2" **\crv{(-5,9) & (5,9)};
  "1";"2" **\crv{(-5,-9) & (5,-9)};
   "1";"2" **\dir{-} ?(.55)*\dir{>};
    (,6.75)*{\scriptstyle >}+(0,3)*{\scriptstyle \phi};
    (0,-6.75)*{\scriptstyle >}+(0,-3)*{\scriptstyle \psi};
    (-10.5,0)*{\scriptstyle \cK};
  (10.5,0)*{\scriptstyle \G};
  {\ar@{=>}_{b}(0,5)*{};(0,1.5)*{}} ;
  {\ar@{=>}_{a}(0,-1.5)*{};(0,-5)*{}} ;
 \endxy.
 \]
}
\newcommand{\hglobe}{
 \[
 \xy
  (-15,0)*{\bullet}="1";
  (0,0)*{\bullet}="2";
  (15,0)*{\bullet}="3";
 "1";"2" **\crv{(-12,9) & (-3,9)};
  "1";"2" **\crv{(-12,-9) & (-3,-9)};
   "2";"3" **\crv{(3,9) & (12,9)} ;
    "2";"3" **\crv{ (3,-9)& (12,-9) };
    (-7.5,6.75)*{\scriptstyle >}+(0,3)*{\scriptstyle \varphi};
    (7.5,6.75)*{\scriptstyle >}+(0,3)*{\scriptstyle \phi};
    (-7.5,-6.75)*{\scriptstyle >}+(0,-3)*{\scriptstyle \varphi'};
    (7.5,-6.75)*{\scriptstyle >}+(0,-3)*{\scriptstyle \phi'};
    (-17.5,0)*{\scriptstyle \cL};
  (17.5,0)*{\scriptstyle \G};
  {\ar@{=>}_{a}(7.5,3)*{};(7.5,-3)*{}} ;
  {\ar@{=>}_{b}(-7.5,3)*{};(-7.5,-3)*{}} ;
   (-18,0)*{};
  (18,0)*{};
 \endxy .
 \]
}
\newcommand{\il}{
\[
\xy
  (-8,0)*{\bullet}="1";
  (8,0)*{\bullet}="2";
 "1";"2" **\crv{(-5,9) & (5,9)};
  "1";"2" **\crv{(-5,-9) & (5,-9)};
   "1";"2" **\dir{-} ?(.55)*\dir{>};
    (,6.75)*{\scriptstyle >}+(0,3)*{\scriptstyle \varphi};
    (0,-6.75)*{\scriptstyle >}+(0,-3)*{\scriptstyle \varphi''};
    (-10.5,0)*{\scriptstyle \cL};
   {\ar@{=>}_{d}(0,5)*{};(0,1.5)*{}} ;
  {\ar@{=>}_{c}(0,-1.5)*{};(0,-5)*{}} ;
     (8,0)*{\bullet}="1";
  (24,0)*{\bullet}="2";
 "1";"2" **\crv{(11,9) & (21,9)};
  "1";"2" **\crv{(11,-9) & (21,-9)};
   "1";"2" **\dir{-} ?(.55)*\dir{>};
    (16,6.75)*{\scriptstyle >}+(0,3)*{\scriptstyle \phi};
    (16,-6.75)*{\scriptstyle >}+(0,-3)*{\scriptstyle \phi''};
  (26.5,0)*{\scriptstyle \G};
  {\ar@{=>}_{b}(16,5)*{};(16,1.5)*{}} ;
  {\ar@{=>}_{a}(16,-1.5)*{};(16,-5)*{}} ;
 \endxy .
 \]
}

\title{On the 1-homotopy type of Lie groupoids}
\author{Hellen Colman}
\address{Department of Mathematics,
Wilbur Wright College,
4300 N. Narragansett Avenue,
Chicago, IL 60634 USA}
\email{hcolman@ccc.edu}
\date{\today}

\begin{document}
\maketitle

\begin{abstract}
We propose a notion of 1-homotopy for generalized maps. This notion generalizes those of natural transformation and ordinary homotopy for functors. The 1-homotopy type of a Lie groupoid is shown to be invariant under Morita equivalence. As an application we consider orbifolds as groupoids and study the notion of orbifold 1-homotopy type  induced by a 1-homotopy between presentations of the orbifold maps.
\end{abstract}

\section*{Introduction}
\label{intro}
We develop a notion of 1-homotopy between generalized maps which is suitable for applications to orbifolds. There are notions of homotopy for functors in general categories that can be applied to groupoids. But they fail to be invariant of Morita equivalence when considering Lie groupoids. Since two Lie groupoids define the same orbifold if they are Morita equivalent,  a notion of homotopy between generalized maps should be invariant of Morita equivalence.

For general categories (no topology nor smooth structure involved) natural transformations play the role of homotopy for functors \cite{Lee,Brown}. Two functors are homotopic if there is a natural transformation between them. We call this notion of homotopy a {\it natural transformation}. 

For topological categories $\cT$ and $\cT'$ the usual notion of homotopy is just a functor which is an ordinary homotopy on objects and on arrows.  We say that two continuous functors $f,g\colon \cT \to \cT'$ are homotopic if there is a continuous functor $H \colon \cT \times \I\to \cT'$ such that $H_0=f$ and $H_1=g$ and $\I$ is the unit groupoid over the interval $[0,1]$ . We call this notion of homotopy an {\it ordinary homotopy}.

Both notions of natural transformation and ordinary homotopy can be adapted to Lie groupoids by requiring all the maps involved to be smooth. None of these two notions is invariant under Morita equivalence.

For an example of Morita equivalent Lie groupoids which are not equivalent by a natural transformation nor an ordinary homotopy, consider the holonomy groupoid $\G=\Hol(M,\F_S)$ associated to a Seifert fibration $\F_S$ on a M\"{o}bius band $M$,  and its reduced holonomy groupoid $\cK=\Hol_T(M,\F_S)$ to a transversal interval $T$. Since the double covering of the M\"{o}bius band by the annulus has no global section, these two groupoids cannot be equivalent by a natural transformation (Section \ref{section}). These groupoids are clearly not equivalent by an ordinary homotopy since the space of objects of $\G=\Hol(M,\F_S)$ is a M\"{o}bius band and the space of objects of $\cK=\Hol_T(M,\F_S)$ is an interval.

We look for a notion of homotopy which is invariant under Morita equivalence and generalizes the notions of natural transformation and ordinary homotopy. To achieve this, we follow a bicategorical approach as emphasized by Landsman in \cite{La1}. Our starting point is the 2-category $\GG$ of Lie groupoids, smooth functors and smooth natural transformations. First we introduce a new notion of 1-homotopy between smooth functors which includes the notions of natural transformation and ordinary homotopy but it is not yet invariant under Morita equivalence. The resulting 2-category is denoted $\HH$. We introduce a notion of {\em essential 1-homotopy equivalence} for the arrows in this 2-category. We prove that the class $W$  of essential 1-homotopy equivalences admits a bicalculus of fractions. The equivalences in this bicategory of fractions $\HH(W^{-1})$ will determine our 1-homotopy equivalences: a generalized map is a {\em 1-homotopy equivalence} if it is an equivalence in $\HH(W^{-1})$.

When inverting the essential equivalences $E$ in the 2-category $\GG$ of Lie groupoids, functors and natural transformations, the following diagram of bicategories commutes:
$$
\xymatrix{ 
\GG
\ar[r]^{U|_{\GG}} 
\ar[d]_{i_{\GG}} 
& \GG(E^{-1})\ar[d]^{i_{\GG(E^{-1})}}\\ 
\HH
\ar[r]^{U} 
& \HH(W^{-1})
}
$$
where $U$ and $U|_{\GG}$  are the universal homomorphisms as defined in \cite{P}. The arrow $i_{\GG(E^{-1})}$ exists by the universal property of $U|_{\GG}$ and can be chosen such that the diagram commutes on the nose.

As a main application we study the 1-homotopy type of orbifolds as groupoids by considering the induced 1-homotopy between orbifold morphisms.

The organization of the paper is as follows. In Section 1 we introduce some basic definitions and constructions for Lie groupoids. In particular, we recall the notion of Morita equivalence and give examples. We introduce some background on bicategories in Section 2 and emphasize the bicategorical viewpoint  by presenting two equivalent constructions of the Morita bicategory $\MM=\GG(E^{-1})$ of Lie groupoids and generalized maps. In Section 3 we present our notion of 1-homotopy for generalized maps. First we recall  the construction of  Haefliger's $\G$-paths and the fundamental groupoid of a groupoid.  Then we introduce the 1-homotopy bicategory $\HH$ of Lie groupoids and strict 1-homotopies, and prove that it admits a right calculus of fractions that inverts the essential 1-homotopy equivalences $W$. Our notion of 1-homotopy between generalized maps corresponds to 2-morphisms in the bicategory $\HH(W^{-1})$. We prove that the notion of 1-homotopy equivalence obtained in this bicategory is invariant under Morita equivalence and generalizes the notions of natural transformation and ordinary homotopy. In Section 4 we study the notion of 1-homotopy for orbifold morphisms induced by a 1-homotopy between their presentation generalized maps.

\subsection*{Acknowledgements}
I am indebted to Peter May and Ieke Moerdijk for helpful suggestions. I also would like to thank Bruce Bartlett for inspiring discussions and Dorette Pronk for her careful reading of an early version of this article and her enlightening comments.

\section{Lie Groupoids}
\subsection{Preliminaries}
A {\it groupoid} $\G$ is a small category in which each arrow is invertible \cite{McL}. Thus a groupoid consists of a set of arrows $G_1$ and a set of objects $G_0$, together with five structure maps:
\begin{enumerate}
\item 
 The maps $\xymatrix{s,t: G_1 \ar@<.5ex>[r] \ar@<-.5ex>[r] & G_0} $
called  source and  target maps. An element $g\in G_1$ with $s(g)=x$ and $t(g)=y$ is an
arrow from $x$ to $y$ and will be denoted by $g:x \rightarrow y$. 

\item The composition map $m: G_1 \times_{G_0} G_1\rightarrow G_1$ mapping a pair $(h,g)$ with $s(h)=t(g)$ to the composition
$m(g,h)=hg$.

\item A unit map $u:G_0\rightarrow G_1$ mapping each $x\in G_0$ to a two-sided unit $u(x)=1_{x}$ for the composition.

\item An inverse map $i: G_1 \rightarrow G_1$ mapping each arrow $g: x \rightarrow y$ in $G$ to a two-sided inverse $i(g)=g^{-1}$ for the composition.
\end{enumerate}
The set of arrows from $x$ to $y$ is denoted $G(x,y)=\{ g\in G_1 | s(g)=x \mbox{ and } t(g)=y\}$. The set of arrows from $x$ to itself is a group called the {\it isotropy} group of $G$ at $x$ and denoted by $G_{x}=G(x,x)$. The {\it orbit} of $x$ is the set $ts^{{-1}}(x)$. The orbit space $|\G|$ of $\G$ is the quotient of $G_{0}$ under the equivalence relation: $x\sim y$ iff $x$ and $y$ are in the same orbit.

A {\it Lie groupoid} $\G$ is a groupoid such that $G_1$ and $G_0$ are manifolds, the structure maps are smooth and $s$ and $t$ are submersions \cite{Mck}. We require the base manifold $G_0$ to be Hausdorff.

\begin{example}
{\it Unit groupoid.} Let $M$ be a smooth manifold. Consider the groupoid $\G$ with $G_0=G_1=M$. This is a Lie groupoid whose arrows are all units, called the unit groupoid and denoted $\G=u(M)$. The orbit space $|\G|$ is again the manifold $M$.
\end{example}

\begin{example}
{\it Pair groupoid.} Let $M$ be a smooth manifold. Consider the groupoid $\G$ with $G_0=M$ and $G_1=M\times M$. This is a Lie groupoid with exactly one arrow from any object $x$ to any object $y$, called the pair groupoid and denoted $\G=\Pair(M)$ or $M\times M$. The orbit space $|\G|$ has only one point.
\end{example}

\begin{example}
{\it Point groupoid.} Let $G$ be a Lie group. Let $\bullet$ be a point. Consider the groupoid $\G$ with $G_0=\bullet$ and $G_1=G$. This is a Lie groupoid with exactly one object $\bullet$ and $G$ is the manifold of arrows in which the maps $s$ and $t$ coincide. We denote the point groupoid by $\bullet^G$. If $G_1$ is also a point, then the point groupoid is called the {\it trivial} groupoid and denoted  $\bf 1$.
\end{example}

\begin{example}
{\it Translation groupoid.} Let $K$ be a Lie group acting (on the left) on a smooth manifold  $M$. Consider the groupoid $\G$ with $G_0=M$ and $G_1=K\times M$. This is a Lie groupoid with arrows $(k, x)$  from any object $x$ to $y=kx$, called the translation or action groupoid and denoted $\G=K\ltimes M$. The orbit space $|\G|$ is the orbit space of the action $M/K$ which is not always a manifold.
\end{example}

\begin{example}
{\it Holonomy groupoid.} Let $(M,\F)$ be a foliated manifold  $M$. Consider the groupoid $\G$ with $G_0=M$ and whose arrows from $x$ to $y$ on the same leaf $L\in\F$ are the holonomy classes of paths in $L$ from $x$ to $y$. There are no arrows between points in different leaves. This is a Lie groupoid \cite{W} called the holonomy groupoid and denoted $\G=\Hol(M,\F)$. The orbit space $|\G|$ is the leaf space of the foliation.
\end{example}

\subsection{Morphisms and equivalences}
From now on all groupoids will be assumed Lie groupoids. The general reference for this section is \cite{M}.

A {\it morphism} $\phi: \cK \to \G$ of groupoids is a functor given by two smooth maps $\phi: K_1 \to G_1$ and $\phi: K_0 \to G_0$ that together commute with all the structure maps of the groupoids $\cK$ and $\G$.

A {\it natural transformation} $T$ between two homomorphisms $\phi, \psi: \cK \to \G$ is a smooth map $T: K_{0} \to G_1$ with $T(x):\phi(x)\to\psi(x)$ such that for any arrow $h:x \rightarrow y$ in $K_1$,  the identity $\psi(h)T(x)=T(y)\phi(h)$ holds. We write $\phi\simT\psi$.

A morphism $\phi: \cK \to \G$ of groupoids is an {\it equivalence} of groupoids if there exists a morphism $\psi: \G \to \cK$ of groupoids and natural transformations $T$ and $T'$ such that $\psi\phi\sim_{T} \id_{\cK}$ and $\phi\psi\sim_{T'} \id_{\G}$.

A morphism $\e: \cK \to \G$ of groupoids is an {\it essential equivalence} of groupoids if
\begin{itemize}
\item[(i)] $\e$ is essentially surjective in the sense that \[t\pi_{1}:G_1\times_{G_0}K_0\rightarrow G_0\] is a surjective submersion
where $G_1\times_{G_0}K_0$ is the pullback along the source $s: G_1\to G_0$;
\item[(ii)] $\e$ is fully faithful in the sense that $K_1$ is  the following pullback of manifolds:
\[\xymatrix{
K_1 \ar[r]^{\e} \ar[d]_{(s,t)}& G_1 \ar[d]^{(s,t)} \\ 
K_0\times K_0 \ar[r]^{\e \times \e} & G_0\times G_0}\]
\end{itemize}
The first condition implies that for any object $y\in G_0$, there exists an object  $x\in K_0$ whose image $\e(x)$ can be connected to $y$ by an arrow $g\in G_1$. The second condition implies that for all $x,z\in K_0$, $\e$ induces a diffeomorphism $K(x,z)\to G(\e(x),\e(z))$ between the submanifolds of arrows.

For general categories the notions of equivalence and essential equivalence coincide. This applies to the particular case in which the categories are groupoids. But when some extra structure is involved (continuity or differentiability) these two notions are not the same anymore. An essential equivalence implies the existence of the inverse functor using the axiom of choice but not the existence of a {\it smooth} functor. The need to invert essential equivalences is what will lead to the definition of generalized maps, a category where essential equivalences can be inverted.

\begin{prop}\cite{MM}
Every  equivalence of Lie groupoids is an essential equivalence.
\end{prop}

\begin{remark}
The converse does not hold for Lie groupoids. 
\end{remark}

{\it Morita equivalence} will be the smallest equivalence relation between Lie groupoids such that they are equivalent whenever there exists an essential equivalence between them. First, we recall the notion of weak pullback. Given the morphisms of groupoids $\psi: \cK \to \G$ and $\phi: \cL \to \G$, the {\it weak pullback} $\cK\times_{\G}\cL$ is a groupoid whose space of objects is 
$$(\cK\times_{\G}\cL)_0=K_0\times_{G_0}G_1\times_{G_0}L_0$$ consisting of triples $(x, g, y)$ with $x\in K_0$, $y\in L_0$ and $g$ an arrow in $G_1$ from $\psi(x)$ to $\phi(y)$. An arrow between $(x, g, y)$ and $(x', g', y')$ is a pair of arrows $(k,l)$ with $k\in K(x,x')$, $l\in L(y,y')$ such that $g'\psi(k)=\phi(l)g$. The space of arrows can be identified with $(\cK\times_{\G}\cL)_1=K_1\times_{G_0}G_1\times_{G_0}L_1$.

If at least one of the two morphisms is a submersion on objects, then the weak pullback $\cK\times_{\G}\cL$ is a Lie groupoid. In this case, the diagram of Lie groupoids
\[\xymatrix{
\cK\times_{\G}\cL \ar[r]^{p_1} \ar[d]_{p_3}& \cK \ar[d]^{\psi} \\ 
\cL \ar[r]^{\phi} & \G}\]
commutes up to a natural transformation and it is universal with this property.
\begin{remark} If $\phi$ is an essential equivalence, then  $\cK\times_{\G}\cL$ is a Lie groupoid and $p_1$ is an essential equivalence too.
\end{remark}

\begin{definition}
Two Lie groupoids $\cK$ and $\G$ are {\it Morita equivalent} if there exists a Lie groupoid $\cJ$ and essential equivalences
\[\cK\overset{\e}{\gets}\cJ\overset{\sigma}{\to}\G.\]
\end{definition}
By means of the weak pullback defined above, we can see that this defines an equivalence relation, denoted $\cK\simM \G$. In this case, it is always possible to chose the equivalences $\e$ and $\sigma$ being surjective submersions on objects \cite{MM}.

\begin{remark}
An essential equivalence of Lie groupoids $\e: \cK \to \G$ induces a homeomorphism $|\e|: | \cK |\to |\G|$ between quotient spaces.
\end{remark}

\begin{example} {\it Reduced holonomy groupoid.} Let $\G=\Hol(M, \F)$ be the holonomy groupoid of a foliation $\F$.
Given a complete transversal $T\subset M$, consider the arrows whose source and target are in $T$. The reduced holonomy groupoid to $T$, $\cK=\Hol_T(M, \F)$,  is defined as the groupoid whose manifold of objects is $K_0=T$ and the manifold of arrows $K_1$ is given by the following pullback of manifolds:

\[\xymatrix{
K_1 \ar[r]^{} \ar[d]_{(s,t)}& G_1 \ar[d]^{(s,t)} \\ 
T\times T\ar[r]^{i} & M\times M}\]
The inclusion functor $\Hol_T(M, \F)\to\Hol(M, \F)$ is an essential equivalence. For a given foliation, all the reduced holonomy groupoids to complete transversals are Morita equivalent.

\end{example}

\section{Groupoid morphisms} The most suitable notion of morphism between Lie groupoids in the context of orbifolds is that of Hilsum-Skandalis maps or generalized maps.    
The composition of these maps is not strictly associative. Groupoids and generalized maps form a bicategory. We review next some of the basic definitions and we set our notation.
\subsection{Bicategories}
A bicategory $\BB$ consists of a class of objects, morphisms between objects and 2-morphisms between morphisms together with various ways of composing them \cite{Be}. We will picture the objects as points:
$$\bullet \scriptstyle \G$$
the morphisms between objects as arrows:
 \[
 \xy
  (-8,0)*{\bullet}="1";
  (8,0)*{\bullet}="2";
 "1";"2" **\dir{-} ?(.55)*\dir{>};
    (0,3)*{\scriptstyle \phi};
    (-10.5,0)*{\scriptstyle \cK};
  (10.5,0)*{\scriptstyle \G};
\endxy
 \]
and the 2-morphisms between morphisms as double arrows:
\globe
\begin{definition} A 2-morphism $a: \phi\ri\psi$ is a {\it 2-isomorphism} in $\BB$ if it is invertible: i.e. if there exists a 2-morphism $b:\psi\ri\phi$ such that $ab=\id_{\psi}$ and $ba=\id_{\phi}$. In this case we will say that the morphisms $\phi$ and $\psi$ are equivalent and write $\phi\sim\psi$.
\end{definition}
\begin{definition} \label{eq} A morphism $\varphi: \cK\to \G$ is an {\it equivalence} in $\BB$ if it is invertible up to a 2-isomorphism: i.e. if there exists a morphism $\xi :\G\to \cK$ such that $\varphi\xi\sim\id_{\G}$ and $\xi\varphi\sim\id_{\cK}$. In this case we will say that the objects $\cK$ and $\G$ are equivalent and write $\cK\sim\G$.
\end{definition}
The composition $\phi\varphi$ of morphisms  $\varphi:\cL\to\cK$ and $\phi:\cK\to\G$ is denoted by:
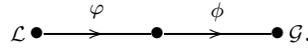
\[
 \xy
  (-8,0)*{\bullet}="1";
  (8,0)*{\bullet}="2";
  (24,0)*{\bullet}="3";
 "1";"2" **\dir{-} ?(.55)*\dir{>};
  "2";"3" **\dir{-} ?(.55)*\dir{>};
    (0,3)*{\scriptstyle \varphi};
     (16,3)*{\scriptstyle \phi};
    (-10.5,0)*{\scriptstyle \cL};
  (26.5,0)*{\scriptstyle \G};
\endxy.
 \]
We can compose 2-morphisms in two ways called horizontal and vertical composition. The horizontal composition $ab$ of 2-morphisms $b:\varphi\ri\varphi'$ and $a:\phi\ri\phi'$ is denoted by:
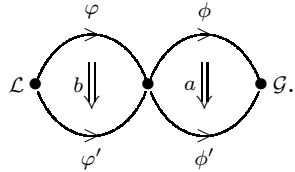
\hglobe
The vertical composition $a\cdot b$ of 2-morphisms  $b:\phi\ri\varphi$ and $a:\varphi\ri\psi$ is denoted by:
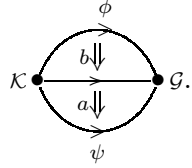
\vglobe
The following interchange law relates horizontal and vertical compositions: $(a\cdot b)(c\cdot d)= (ac)\cdot (bd)$
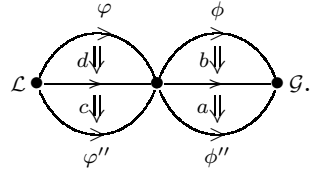
\il
Vertical composition is strictly associative whereas horizontal composition is only associative up to 2-isomorphism (associator). The unit laws for morphisms hold up to 2-isomorphisms (left and right unit constraints). Associator and unit constrains are required to be natural with respect to their arguments and verify certain coherence axioms \cite{Be}.

A 2-category is a bicategory in which the natural 2-isomorphisms above are identities. The category of Lie groupoids and functors can be seen as a 2-category $\GG$ with natural transformations as 2-morphisms. We will show two different constructions producing equivalent  bicategories $\MM$ and $\MM'$ in which the morphisms are respectively generalized maps and bibundles.

\subsection{Generalized maps}
We describe in this section the generalized maps obtained by localization of essential equivalences  \cite{P,GZ}. Considering the bicategory $\GG$ of Lie groupoids, functors and natural transformations, the bicategory $\MM$ is obtained as the bicategory of fractions of $\GG$ when inverting the essential equivalences $E$, $\MM=\GG(E^{-1})$.

A {\it generalized map} from $\cK$ to $\G$ is a pair of morphisms 
\[\cK\overset{\e}{\gets}\cJ\overset{\phi}{\to}\G\]
such that $\e$ is an essential equivalence. We denote a generalized map by $(\e,\phi)$. Roughly speaking, a generalized map  from $\cK$ to $\G$ is obtained by first replacing $\cK$ by another groupoid $\cJ$ essentially equivalent to it and then mapping $\cJ$ into $\G$ by an ordinary morphism.

Two generalized maps from $\cK$ to $\G$, 
$\cK\overset{\e}{\gets}\cJ\overset{\phi}{\to}\G$ and 
$\cK\overset{\e'}{\gets}\cJ'\overset{\phi'}{\to}\G$, 
are {\it isomorphic} if there exists a groupoid $\cL$ and essential equivalences \[\cJ\overset{\alpha}{\gets}\cL\overset{\beta}{\to}\cJ'\] such that the diagram

$$
\xymatrix{ & 
{\cJ}\ar[dr]^{\phi}="0" \ar[dl]_{\e}="2"&\\
{\cK}&{\cL} \ar[u]_{\alpha} \ar[d]^{\beta}&{\G}\\
&{\cJ'}\ar[ul]^{{\e}'}="3" \ar[ur]_{\phi'}="1"&
\ar@{}"0";"1"|(.4){\,}="7" 
\ar@{}"0";"1"|(.6){\,}="8" 
\ar@{}"7" ;"8"_{\sim_{T'}} 
\ar@{}"2";"3"|(.4){\,}="5" 
\ar@{}"2";"3"|(.6){\,}="6" 
\ar@{}"5" ;"6"^{\sim_T} 
}
$$
commutes up to natural transformations. We write $(\e,\phi)\sim(\e',\phi')$.

In other words, there are natural transformations $T$ and $T'$ such that the generalized maps
\[\cK\overset{\e\alpha}{\gets}\cL\overset{\phi\alpha}{\to}\G\]
\[\cK\overset{\e'\beta}{\gets}\cL\overset{\phi'\beta}{\to}\G\] satisfy $\e\alpha\simT\e'\beta$ and $\phi\alpha\sim_{T'}\phi'\beta$.

\begin{remark}\label{remark}
\begin{enumerate}
\item If $\cK\overset{\e}{\gets}\cJ\overset{\phi}{\to}\G$ and 
$\cK\overset{\e}{\gets}\cJ\overset{\phi'}{\to}\G$
are two generalized maps with  $\phi\sim_{T}\phi'$ then $(\e,\phi)\sim(\e,\phi')$.
\item If $\cK\overset{\e}{\gets}\cJ\overset{\phi}{\to}\G$ and 
$\cK\overset{\e'}{\gets}\cJ'\overset{\phi'}{\to}\G$
are two generalized maps and $\delta:\cJ'\to \cJ$ an essential equivalence with $\phi'=\phi\delta$ and $\e'=\e\delta$ then $(\e,\phi)\sim(\e',\phi')$.
\end{enumerate}
\end{remark}
There is an equivalence relation between the diagrams above (see \cite{Pthesis}). A {\em 2-isomorphism} is an equivalence class of diagrams. Vertical and horizontal composition of diagrams are defined as in \cite{Pthesis}.

\begin{prop} \cite{P} The collection of Lie groupoids as objects, generalized maps as morphisms and  2-isomorphisms is a  bicategory.
\end{prop}
This bicategory will be denoted by $\MM$. All the 2-morphisms in $\MM$ are isomorphisms.

For each groupoid $\G$ the unit arrow $(\id, \id)$ is defined as  the generalized map 
$\displaystyle\G\xleftarrow{\id}\G\xrightarrow{\id}\G.$
The composition of two arrows $(\G\overset{\delta}{\gets}\cJ'\overset{\varphi}{\to}\cL)\circ (\cK\overset{\e}{\gets}\cJ\overset{\phi}{\to}\G)$
 is given by  the generalized map:
\[\cK\xleftarrow{\e p_1}\cJ\times_{\G}\cJ'\xrightarrow{\varphi p_3}\cL\]
where $p_1$ and $p_3$ are the projections in the following weak pullback of groupoids:
\[\xymatrix{
\cJ\times_{\G}\cJ' \ar[r]^{p_1} \ar[d]_{p_3}& \cJ \ar[d]^{\phi} \ar[r]^{\e}&\cK \\ 
\cJ' \ar[r]^{\delta}\ar[d]^{\varphi} & \G\\
\cL}\]
The morphism $p_1$ is an essential equivalence since it is the weak pullback of the essential equivalence $\delta$. Then $\e p_1$ is an essential equivalence. This composition is associative up to isomorphism.

The unit arrow is a left and right unit for this multiplication of arrows up to isomorphism. The composition 
$(\G\xleftarrow{\id}\G\xrightarrow{\id}\G)\circ (\cK\overset{\e}{\gets}\cJ\overset{\phi}{\to}\G)$ is the generalized map 
$\cK\xleftarrow{\e p_1}\cJ\times_{\G}\G\xrightarrow{\varphi p_3}\G$. Since $\varphi=\delta=\id$ implies  that $\varphi p_3=\phi p_1$ and $p_1$ is an essential equivalence. We have that 
$(\e p_1,\phi p_1)\sim(\e,\phi)$  by Remark \ref{remark} (2). 

We will present next another description of these generalized maps in terms of groupoid bundles which provides a more concrete approach.

\subsection{Hilsum-Skandalis maps}
These maps were introduced by Haefliger in \cite{H1} and developed further in different contexts  \cite{HS}, \cite{Pra}, \cite{M1}. Mr\v{c}un studied these maps for general groupoids in his 1996 thesis \cite{Mr}, on which work this section is largely based.
\subsubsection{Actions of groupoids on manifolds}
Let $M$ be a manifold, $\G$ a groupoid and $\mu:M\to G_0$ a smooth map. 
A {\it right action} of $\G$ on $M$ is a map 
$$M\times_{G_0}^tG_1{\to}M,\qquad (x,g)\mapsto x g$$
defined on $M\times_{G_0}^tG_1$ given by the following pullback of manifolds along the target map:
\[\xymatrix{
M\times_{G_0}^tG_1 \ar[r]^{p_1} \ar[d]_{p_2}& M \ar[d]^{\mu} \\ 
G_1 \ar[r]^{t} & G_0}\]
such that $\mu(xg)=t(g)$, $x1=x$, $(xg)h=x(gh)$.

Analogously, we have a {\it left action} by considering the pullback $G_1\times_{G_0}^sM$ along the source map.

The {\it translation groupoid} $M\rtimes\G$ associated to a right action of $\G$ on $M$ is given by
$(M\rtimes\G)_0=M$ and 
$(M\rtimes\G)_1=M\times_{G_0}^tG_1$ where the source map is given by the action $s(x, g)=xg$ and the target map is just the projection $t(x,g)=x$.

The {\it double translation groupoid} $\cK\ltimes M\rtimes\G$ associated to a left action of $\cK$ on $M$ and a right action of $\G$ on $M$ which commute with each other is given by
$(\cK\ltimes M\rtimes\G)_0=M$ and 
$(\cK\ltimes M\rtimes\G)_1=K_1\times_{K_0}^s M\times_{G_0}^tG_1$ where the space of arrows is obtained by the following pullbacks of manifolds:

\[\xymatrix{
K_1\times_{K_0}^s M\times_{G_0}^tG_1 \ar[rr]^{}\ar[d]^{} &&K_1\ar[d]^{s}\\
M\times_{G_0}^tG_1 \ar[r]^{p_1} \ar[d]_{p_2}& M \ar[d]^{\rho} \ar[r]^{\tau}&K_0 \\ 
G_1\ar[r]^{t}& G_0}\] 
then $K_1\times_{K_0}^s M\times_{G_0}^tG_1 =\{(h,x,g)\; |\; s(h)=\tau(x) \mbox{ and } t(g)=\rho(x)\}$ with $s(h,x,g)=x$ and  $t(h,x,g)=hxg^{-1}$. The composition of arrows is given by $(h,x,g)(h',x',g')=(hh',x',gg')$.

\subsubsection{Bibundles}
A {\it right $\G$-bundle} $M$ over $B$ is a map $\pi: M\to B$ with a right action of $\G$ on $M$ preserving the fibers. A right $\G$-bundle is {\it principal} if the map $$M\times_{G_0}^tG_1 \overset{\alpha}{\to}M\times_{B}M,\qquad \alpha(x,g)=(x g, x)$$ is a diffeomorphism and $\pi$ is a surjective submersion.

A {\it $\cK\G$-bibundle} $M$ is a left $\cK$-bundle over $G_0$ as well as a right $\G$-bundle over $K_0$. We represent  a $\cK\G$-bibundle by the following diagram:
$$
\begin{diagram}
\node[2]{M}
        \arrow{se,l}{\tau} \arrow{sw,l}{\rho}\\
\node{G_0}\node[2]{K_0}
\end{diagram}
$$
where $\rho$ is a left $\cK$-bundle and $\tau$ is a right $\G$-bundle. We denote a $\cK\G$-bibundle $M$ as $(\cK, M, \G)$ and we write  $(\cK, \rho, M, \G, \tau)$ if we need to specify the bundle maps.

A $\cK\G$-bibundle $M$ is {\it right principal} if the right $\G$-bundle $\tau: M\to K_0$ is principal. In this case,  $M\times_{G_0}^tG_1$ is diffeomorphic to 
$M\times_{K_0}M$ and $M/\G$ is diffeomorphic to $K_0$. Analogously, left principal.

A $\cK\G$-bibundle $M$ is {\it biprincipal} if both left and right bundles $\rho$ and $\tau$ are principal.

Two $\cK\G$-bibundles $M$ and $N$ are {\it isomorphic} if there is a diffeomorphism $f:M\to N$ that intertwines the maps $M\to G_0$, $M\to K_0$ with the maps $N\to G_0$, $N\to K_0$ and also intertwines the $\cK$ and $\G$ actions. In other words,  $f(hxg)=hf(x)g$ and $\tau=\tau'f$, $\rho=\rho'f$. We write $(\cK, M, \G)\sim(\cK, N, \G)$.
$$
\begin{diagram}
\node[2]{M}
        \arrow{se,l}{\tau} \arrow{sw,l}{\rho}\arrow[2]{s,l}{f}\\
\node{G_0}\node[2]{K_0}\\
\node[2]{N}
        \arrow{ne,r}{\tau'} \arrow{nw,r}{\rho'}\\
\end{diagram}
$$

\begin{definition} A Hilsum-Skandalis map $|(\cK, M, \G)|$ is an isomorphism class of right principal $\cK\G$-bibundles.
\end{definition}

\begin{prop}\cite{La1} The collection of all Lie groupoids as objects, right principal $\cK\G$-bibundles as morphisms and isomorphisms $f$ as 2-morphisms is a bicategory. 
\end{prop}
This bicategory will be denoted $\MM'$. All 2-morphisms in $\MM'$ are isomorphisms.
For each groupoid $\G$ the unit arrow is defined as the $\G\G$-bibundle
$$
\begin{diagram}
\node[2]{G_1}
        \arrow{se,l}{s} \arrow{sw,l}{t}\\
\node{G_0}\node[2]{G_0}
\end{diagram}
$$
The left and right actions of $\G$ on $G_1$ are given by the multiplication in the groupoid $\G$.

The multiplication of arrows $(\cK, M, \G)$ and $(\G, N, \cL)$ is given by  the bibundle $(\cK, (M\times_{G_0}N)/\G, \cL)$ where $M\times_{G_0}N$ is the pullback of manifolds
\[\xymatrix{
M\times_{G_0}N \ar[r]^{  \;p_1} \ar[d]_{p_2}& M \ar[d]^{\rho_M} \\ 
N \ar[r]^{\tau_N} & G_0}\]
and in addition, $\G$ acts on the manifold $M\times_{G_0}N$ on the right by $(x,y) g=(xg,g^{-1}y)$. The orbit space is a $\cK\cL$-bibundle
$$
\begin{diagram}
\node[2]{(M\times_{G_0}N)/\G}
        \arrow{se,l}{\tau} \arrow{sw,l}{\rho}\\
\node{L_0}\node[2]{K_0}
\end{diagram}
$$
where $\tau([x,y])=\tau_M(x)$ and $\rho([x,y])=\rho_N(y)$. The left $\cK$-action is given by $k[x,y]=[kx,y]$ and the right $\cL$-action by $[x,y]l=[x,yl]$. This bibundle is right principal. The multiplication is associative up to isomorphism.

The unit arrow   $(\G, G_1, \G)$   is a left and right unit for this multiplication of arrows up to isomorphism. We have that the bibundle $(\cK, (M\times_{G_0}^tG_1)/\G, \G)$ is isomorphic to $(\cK, M, \G)$ since  the map 
$$(M\times_{G_0}^tG_1)/\G \overset{f}{\to}M,\qquad f([x,y])=xy$$ is a diffeomorphism satisfying $f(h[x,y]g)=hf([x,y])g$. Hence there is a 2-morphism $f$ from $(\cK, M, \G)$ to the composition $(\G, G_1, \G)\circ (\cK, M, \G)= (\cK, (M\times_{G_0}^tG_1)/\G, \G)$.

\subsection{A biequivalence $\Gamma\colon \MM'\to\MM$ of bicategories}
We will show an explicit construction of a bijective correspondence between generalized maps and bibundles which will allow us to switch from one formulation to the other when needed. This extends to generalized maps the construction given in \cite{La1}, \cite{La} and \cite{MM1} for strict maps.

In addition $\MM'$ is biequivalent to $\MM$. Recall that a homomorphism of bicategories is a generalization of the notion of a functor sending objects, morphisms and 2-morphisms of one bicategory to items of the same types in the other one, preserving compositions and units up to 2-isomorphism \cite{Be}. A homomorphism  $\Gamma\colon \MM'\to\MM$ is a biequivalence if the functors $\MM'(\cK,\G)\to  \MM(\Gamma\cK,\Gamma\G)$ are equivalences for all objects $\cK$ and $\G$ of $\MM'$ and if for every object $\cL$ of $\MM$ there is an object $\cK$ of $\MM'$ such that $\Gamma\cK$ is equivalent to $\cL$ in $\MM$.

Given a right principal $\cK\G$-bibundle $M$:
$$
\begin{diagram}
\node[2]{M}
        \arrow{se,l}{\tau} \arrow{sw,l}{\rho}\\
\node{G_0}\node[2]{K_0}
\end{diagram}
$$
where $\rho$ is a left $\cK$-bundle and $\tau$ is a right principal $\G$-bundle we construct a generalized map 
\[\cK\overset{\e}{\gets}\cJ\overset{\phi}{\to}\G\]
by taking $\cJ=\cK\ltimes M\rtimes\G$ and the following morphisms $\e$ and $\phi$:
$$M\overset{\e_0}{\to}K_{0},\qquad \e_0=\tau \;\; \mbox{  and  } \;\;
K_1\times_{K_0}^sM\times_{G_0}^tG_1 \overset{\e_1}{\to}K_{1},\qquad \e_1=p_1$$
$$M\overset{\phi_0}{\to}G_{0},\qquad \phi_0=\rho\;\; \mbox{  and  } \;\; K_1\times_{K_0}^sM\times_{G_0}^tG_1 \overset{\phi_1}{\to}G_{1},\qquad \phi_1=p_3$$
since $\tau$ is a principal bundle, we have that $\e$ is an essential equivalence.

We will show that if $(\cK, M, \G)\sim(\cK, N, \G)$ then the associated generalized maps $(\e,\phi)$ and $(\e',\phi')$ are isomorphic. Let $f: M\to N$ be the equivariant diffeomorphism that intertwines the bundles. Define $$\bar f :  \cK\ltimes M\rtimes\G{\to}\cK\ltimes N\rtimes\G$$ by $\bar f_0=f$ on objects and $\bar f_1(h,x,g)=h x g$ on arrows. These maps commute with all the structural maps by the equivariance of  $f$. Since $\bar f_0$ is a diffeomorphism, it is in particular a surjective submersion and the manifold of arrows $K_1\times_{K_0}^sM\times_{G_0}^tG_1$ is obtained from the following pullback of manifolds:
\[\xymatrix{
K_1\times_{K_0}^sM\times_{G_0}^tG_1 \ar[r]^{\bar f_1} \ar[d]_{(s,t)}&K_1\times_{K_0}^sN\times_{G_0}^tG_1  \ar[d]^{(s,t)} \\ 
M\times M \ar[r]^{\bar f_0\times\bar f_0}&N\times N}\]
Then $\bar f$ is an essential equivalence. Also, as $f$ intertwines the bundles, we have that  $\phi'=\phi\bar f$ and $\e'=\e\bar f$ and by Remark \ref{remark} (2)  follows that $(\e,\phi)\sim(\e',\phi')$.

We define a homomorphism of bicategories by
$$\Gamma: \MM' \to\MM ,\qquad \Gamma ((\cK, M, \G))=(\e,\phi)$$ as constructed above on morphisms and being the identity map  on objects.

For 2-morphisms $f: M\to N$, we define $\Gamma(f)$ as the following diagram:
$$
\begin{diagram}
\node[2]{ \cK\ltimes M\rtimes\G}
        \arrow{se,l}{\phi} \arrow{sw,l}{\e}\\
\node{\cK}\node{\cK\ltimes M\rtimes\G} \arrow{n,r}{\id} \arrow{s,r}{\beta}\node{\G}\\
\node[2]{\cK\ltimes N\rtimes\G}\arrow{nw,r}{\e'} \arrow{ne,r}{\phi'}
\end{diagram}
$$
where $\beta: \cK\ltimes M\rtimes\G \to \cK\ltimes N\rtimes\G$ is defined by $\beta(x)=f(x)$ on objects and $\beta(h,x,g)=(h,f(x),g)$ on arrows. Since $\tau'f=\tau$ and $\rho'f=\rho$ we have that 
$s(h)=\tau'(f(x))$ and $t(g)=\rho'(f(x))$.

Conversely, given a generalized map  from $\cK$ to $\G$
 \[\cK\overset{\e}{\gets}\cJ\overset{\phi}{\to}\G\]
we construct an associated right principal  $\cK\G$-bibundle $M$
$$
\begin{diagram}
\node[2]{M}
        \arrow{se,l}{\tau} \arrow{sw,l}{\rho}\\
\node{G_0}\node[2]{K_0}
\end{diagram}
$$
where $M$ is  the quotient by the action of $\cJ$ on $\ti M=J_0\times_{G_0}^tG_1\times_{K_0}^tK_1$  given by the following pullbacks of manifolds:

\[\xymatrix{
\ti M\ar@/_2.3pc/[dd]_{p_4}  \ar@/^1.5pc/[rr]^{p_2}\ar[r]^{}\ar[d]^{} &J_0\times_{G_0}G_1 \ar[r]^{}\ar[d]^{}&G_1 \ar[d]^{t}\ar[r]^{s}&G_0\\
J_0\times_{K_0}K_1\ar[r]^{} \ar[d]_{}&J_0 \ar[r]^{\e} \ar[d]^{\phi}& G_0  &\\ 
K_1\ar[r]^{t} \ar[d]^{s}&K_0&&\\
K_0&&&
}\]
 The maps $\rho$ and $\tau$ are induced in the quotient by $\ti \rho=s p_4$ and $\ti \tau=sp_2$. The action of $\cJ$ on $\ti M$ is given by
 $ ((a,b,d), j)\mapsto (t(j),b\phi(j),d\e(j)).$
The left action of $\cK$ on $M=\ti M/\cJ$ is given by
$$K_1\times_{K_0}^sM\to M ,\qquad (h,[a,b,d])\mapsto [a, bk^{-1}, c,d]$$ and the right action of $\G$ by 
$$M\times_{G_0}^t G_1\to M ,\qquad ([a,b,d]), g)\mapsto [a,b,dg].$$
If $(\e,\phi)\sim(\e',\phi')$ then the associated bibundles  $(\cK, M, \G)$ and $(\cK, N, \G)$ are isomorphic.

The homomorphism $\Gamma$ preserves compositions and units up to 2-isomorphism and it is a biequivalence.

\subsection{Strict maps} We will use now this explicit construction to characterize the generalized maps that come from a strict map.

Any strict morphism $\phi:\cK \to \G$ can be viewed as a generalized map by $$\cK\xleftarrow{\id}\cK\xrightarrow{\phi}\G.$$
The corresponding bibundle is constructed by taking $$M=(K_0\times_{G_0}^t G_1\times_{K_0}^t K_1)/\cK.$$ 
Since
$$(K_0\times_{G_0}^t G_1\times_{K_0}^t K_1)/\cK=K_0\times_{G_0}^tG_1$$ 
$$[a,b,d]\mapsto (t(d),b\phi(d))$$
this bibundle is isomorphic to 
$$
\begin{diagram}
\node[2]{K_0\times_{G_0}^t G_1}
        \arrow{se,l}{p_1} \arrow{sw,l}{sp_2}\\
\node{G_0}\node[2]{K_0}
\end{diagram}
$$

\begin{prop}\cite{Mr} Let $\cK\xleftarrow{\e}\cJ\xrightarrow{\phi}\G$ be a generalized map  and 
$(\cK, \rho, M, \G, \tau)$
 the associated right principal bibundle. Then $\phi$ is an essential equivalence iff $\rho$ is principal.
\end{prop}
\begin{remark} Essential equivalences correspond to {\it biprincipal} bibundles.
\end{remark}
In both bicategories $\MM$ and $\MM'$, Morita equivalences are the invertible morphisms up to a 2-isomorphism, i.e. the {\it equivalences} in $\MM$ (or $\MM'$). If $\cK\simM\G$, let 
$\cK\xleftarrow{\e}\cJ\xrightarrow{\delta}\G$ be the associated generalized map  in $\MM$ with $\e$ and $\delta$ essential equivalences, then the inverse generalized map  is
$\G\xleftarrow{\delta}\cJ\xrightarrow{\e}\cK$. 

In the category $\MM'$, let
$$
\begin{diagram}
\node[2]{M}
        \arrow{se,l}{\tau} \arrow{sw,l}{\rho}\\
\node{G_0}\node[2]{K_0}
\end{diagram}
$$ be the biprincipal $\cK\G$-bibundle representing the Morita equivalence, then the inverse biprincipal $\G\cK$-bibundle is 
$$
\begin{diagram}
\node[2]{M}
        \arrow{se,l}{\rho} \arrow{sw,l}{\tau}\\
\node{K_0}\node[2]{G_0}
\end{diagram}
$$ where the new actions are obtained from the original ones composing with the inverse: the left action of $\G$ on $M$ is given by $g\ast x=xg^{-1}$ induced by the right action of $\G$ on the original bibundle. Similarly, the left action of $\cK$ induces a right action in the inverse bundle.

\begin{prop} \cite{Mr} Let $(\cK, \rho, M, \G, \tau)$ be a right principal $\cK\G$-bibundle and $(\e,\varphi)=\Gamma ((\cK, M, \G))$ its associated generalized map. Then 
$(\e,\varphi)\sim(\id,\phi)$ if and only if $\tau$ has a section.
\end{prop}
In other words, a generalized map comes from a strict map iff when seen as a bibundle, the right principal $\G$-bundle has a section.

If a strict map $\e: \cK\to \G$ is an essential equivalence, regarded as a generalized map it will be invertible. The inverse of a generalized (strict) map $$\cK\overset{\id}{\gets}\cJ\overset{\e}{\to}\G$$ is the generalized map $\G\overset{\e}{\gets}\cJ\overset{\id}{\to}\cK$ which will not always come from a strict map.

\begin{example}\label{section}
Consider the holonomy groupoid associated to the Seifert fibration on the M\"{o}bius band, $\G=\Hol(M,\F_S)$, and the reduced holonomy groupoid $\cK=\Hol_T(M,\F_S)$ to a transversal $T$ given by an interval $I$ transversal to the leaves. The inclusion functor $i_{\cK}: \cK \hookrightarrow  \G$ is an essential equivalence. The associated $\cK\G$-right principal bibundle is 
\[\xymatrix{
&I\times_M^t(M\times S^1){=I\times S^1}\ar[ld]_{\rho}\ar[rd]^{\tau}&  \\ 
M & &I {\ar@/_2pc/[lu]}}\]   
which has a section since it comes from the strict map $i_{\cK}$.

The inverse $\G\cK$-bibundle is 
\[\xymatrix{
&I\times_M^t(M\times S^1){=I\times S^1}\ar[ld]_{\tau}\ar[rd]^{\rho}&  \\ 
I & &M }\]   
where the map $\rho: I\times S^1\to M$ is the double covering of the M\"{o}bius band by an annulus,  which does not have a section. So the corresponding generalized map is not isomorphic to a strict map.
This example shows that the strict map $i_{\cK}$ has an inverse which is not an strict map. The groupoids $\cK$ and $\G$ are Morita equivalent but they are not 
equivalent by a natural transformation.
\end{example}

\section{A 1-homotopy notion}
In this section, we introduce a notion of {\it 1-homotopy equivalence} and {\it essential 1-homotopy equivalence}  for strict maps. Then we prove that the essential 1-homotopy equivalences admit a bicalculus of fractions, and we define our notion of 1-homotopy between generalized maps as a 2-arrow in the bicategory of fractions.

\subsection{Haefliger paths}
We first recall the notions of equivalence and homotopy of $\G$-paths due to Haefliger \cite{H2,H3}.
Let $\G$ be a Lie groupoid and $x$, $y$ objects in $G_0$. A {\it $\G$-path} from $x$ to $y$ over a subdivision $0=t_0\le t_1\le\cdots\le t_n=1$ is a sequence:

$$(g_0,\al_1,g_1,\ldots, \al_n, g_n)$$
where 
\begin{enumerate}
\item $\al_i:[t_{i-1}, t_i]\to G_0$ is a path for all $1\le i\le n$ and
\item $g_i\in G_1$ is an arrow such that 

 $s(g_0)=x$ and $t(g_n)=y$ 

$s(g_i)=\al_i(t_i)$  for all $0< i\le n$

$t(g_i)=\al_{i+1}(t_i)$ for all $0\le i< n$

\[
\xy
 (-7,0)*{\bullet}="1";
 (7,0)*{\bullet}="2";
  (21,0)*{\bullet}="3";
   (35,0)*{\bullet}="4";
 (49,0)*{\bullet}="5";
  (63,0)*{\bullet}="6";
  (77,0)*{\bullet}="7";
   (-10,0)*{x};
    (80,0)*{y};
  (0,2.5)*{\scriptstyle g_0};
  (14,2.5)*{\scriptstyle \al_1};
   (28,2.5)*{\scriptstyle g_1};
   (56,2.5)*{\scriptstyle \al_n};
    (70,2.5)*{\scriptstyle g_n};
  "1"; "2" **\dir{-} ?(.55)*\dir{>};
   "3"; "4" **\dir{-} ?(.55)*\dir{>};
    "6"; "7" **\dir{-} ?(.55)*\dir{>};
 { \ar@{~>}"2"; "3"};
 { \ar@{~>}"5"; "6"};
{ \ar@{}"4"; "5"};
(42,0)*{\cdots};
\endxy 
\]
\end{enumerate}
We define an {\it equivalence} relation $\sim$ among $\G$-paths generated by the following operations:
\begin{enumerate}
\item Add a new point $s\in[t_{i-1}, t_i]$ to the subdivision, take the restrictions  $\al'_i$ and $\al''_i$ of the corresponding path $\al_i$ to the new intervals $[t_{i-1}, s]$ and $[s, t_i]$ and add the identity arrow $1_{\al(s)}$

\medskip
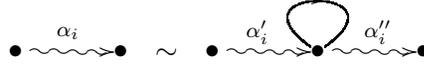
\begin{center}
$\xy
 (-7,0)*{\bullet}="1";
 (7,0)*{\bullet}="2";
    (0,2.5)*{\scriptstyle \al_i};
  { \ar@{~>}"1"; "2"};
  \endxy 
  \:\:\:\sim\:\:\:
   \xy
 (14,0)*{\bullet}="1";
 (28,0)*{\bullet}="2";
  (42,0)*{\bullet}="3";
  (20,2.6)*{\scriptstyle \al'_i};
  (36,2.6)*{\scriptstyle \al''_i};
      "2";"2" **\crv{(14,9)& (42,9)}?(.55)*\dir{>};
   { \ar@{~>}"1"; "2"};
   { \ar@{~>}"2"; "3"};
\endxy $
\end{center}

\item Given a map $h_i:[t_{i-1}, t_i]\to G_1$ with $s\circ h_i=\al_i$, replace: 

$\al_i$ by $t\circ h_i$

$g_{i-1}$ by $h_i(t_{i-1})g_{i-1}$  and 

$g_{i}$ by $g_i (h_i(t_{i}))^{-1}$

\medskip
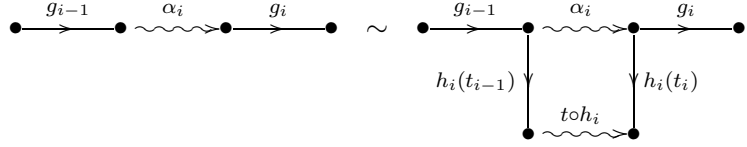
\begin{center}
$\xy
 (13,0)*{\bullet}="1";
 (27,0)*{\bullet}="2";
  (41,0)*{\bullet}="3";
   (55,0)*{\bullet}="4";
  (20,2.5)*{\scriptstyle g_{i-1}};
  (34,2.5)*{{\scriptstyle \al_i}};
   (48,2.5)*{{\scriptstyle g_i}};
    "1"; "2" **\dir{-} ?(.55)*\dir{>};
   "3"; "4" **\dir{-} ?(.55)*\dir{>};
   { \ar@{~>}"2"; "3"};
 \endxy 
  \:\:\:\sim\:\:\:
   \xy
 (13,0)*{\bullet}="1";
 (27,0)*{\bullet}="2";
  (41,0)*{\bullet}="3";
   (55,0)*{\bullet}="4";
   (27,-14)*{\bullet}="5";
   (41,-14)*{\bullet}="6";
 (20,2.5)*{\scriptstyle g_{i-1}};
  (34,2.5)*{{\scriptstyle \al_i}};
   (48,2.5)*{{\scriptstyle g_i}};
   (34,-11.5)*{\scriptstyle t\circ h_{i}};
  (20,-7)*{\scriptstyle h_i(t_{i-1})};
   (46,-7)*{\scriptstyle h_i(t_{i})};
    "1"; "2" **\dir{-} ?(.55)*\dir{>};
   "3"; "4" **\dir{-} ?(.55)*\dir{>};
    "2"; "5" **\dir{-} ?(.55)*\dir{>};
     "3"; "6" **\dir{-} ?(.55)*\dir{>};
   { \ar@{~>}"2"; "3"};
   { \ar@{~>}"5"; "6"};

\endxy $
\end{center}
\end{enumerate}

\begin{remark}
Note that equivalence classes of $\G$-paths correspond to isomorphism classes of generalized maps from $\I$ to $\G$, where $\I$ is the unit groupoid associated to the interval $I=[0,1]$.
\end{remark}

\bigskip
A {\it deformation} between the $\G$-paths  $(g_0,\al_1,g_1,\ldots, \al_n, g_n)$ and  $(g'_0,\al'_1,g'_1,\ldots, \al'_n, g'_n)$ from $x$ to $y$ is given by homotopies
$$H_i:[t_{i-1}, t_i]\times I \to G_0 \mbox { with }(H_i)_0=\al_i \mbox { and }(H_i)_1=\al'_i$$ for $i=1,\ldots, n$ and 
$$\gamma_i: I \to G_1   \mbox { with }(\gamma_i)_0=g_i \mbox { and }(\gamma_i)_1=g'_i$$
for $i=1,\ldots, n-1$, such that  $(g_0,(H_1)_s,(\gamma_1)_s,\ldots,(\gamma_{n-1})_s, (H_n)_s, g_n)$ is a $\G$-path for each $s\in I$
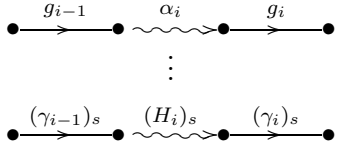
\[
 \xy
 (13,0)*{\bullet}="1";
 (27,0)*{\bullet}="2";
  (41,0)*{\bullet}="3";
   (55,0)*{\bullet}="4";
   (27,-14)*{\bullet}="5";
   (41,-14)*{\bullet}="6";
   (13,-14)*{\bullet}="7";
   (55,-14)*{\bullet}="8";
 (20,2.5)*{\scriptstyle g_{i-1}};
  (34,2.5)*{{\scriptstyle \al_i}};
   (34,-4)*{{\vdots}};
   (48,2.5)*{{\scriptstyle g_i}};
    (20,-11.5)*{\scriptstyle (\gamma_{i-1})_s};
  (34,-11.5)*{{\scriptstyle (H_i)_s}};
   (48,-11.5)*{{\scriptstyle  (\gamma_{i})_s}};
    "1"; "2" **\dir{-} ?(.55)*\dir{>};
   "3"; "4" **\dir{-} ?(.55)*\dir{>};
      "7"; "5" **\dir{-} ?(.55)*\dir{>};
     "6"; "8" **\dir{-} ?(.55)*\dir{>};
   { \ar@{~>}"2"; "3"};
   { \ar@{~>}"5"; "6"};
\endxy
 \]

\begin{definition} Two $\G$-paths between $x$ and $y$ are {\it homotopic} if one can be obtained from the other by a sequence of equivalences and deformations.
\end{definition}
We define a multiplication of homotopy classes of  $\G$-paths by $$[(g'_0,\al'_1,\ldots, \al'_n, g'_n)]\;[(g_0,\al_1,\ldots, \al_n, g_n)]=[g_0,\al_1,\ldots,\al_n,g'_0 g_n,\al'_1,\ldots,\al'_n,g'_n]$$
where $g'_0g_n$ is the multiplication of two composable arrows in $\G$ and the paths $\al_i$ are reparametrized to the new subdivision.

The inverse of the homotopy class of the $\G$-path $[(g_0,\al_1,\ldots, \al_n, g_n)]$ from $x$ to $y$ is the class of the $\G$-path from $y$ to $x$
$$[(g\inv_k,\al'_1,\ldots, g\inv_1,\al'_k, g\inv_0)]$$ over the same subdivision where $\al'_i:[t_{i-1},t_i]\to G_0$ is given by 
$$\al'_i(t)=\al_{k-i+1}\left(t_{k-i+1}+\left({t_{k-i}-t_{k-i+1}\over t_{i-1}-t_i}\right)(t_{i-1}-t)\right)$$

\begin{definition}\cite{MM1}  The fundamental groupoid $\pi_1(\G)$ of the Lie groupoid $\G$ is a groupoid over  $G_0$ whose arrows are the homotopy classes of 
$\G$-paths with the multiplication defined above.
\end{definition}

We will denote the fundamental groupoid $\pi_1(\G)$ also by $\G_*$.

\begin{prop} \cite{MM1} The fundamental groupoid $\G_*$ of a Lie groupoid $\G$ is a Lie groupoid.
\end{prop}

A morphism  $\phi: \cK\to \G$ of Lie groupoids induces a morphism $\phi_*: \cK_*\to \G_*$ between the fundamental groupoids given by $\phi_*=\phi$ on objects and 
$$\phi_*([g_0,\al_1,g_1,\ldots, \al_n, g_n])=[\phi(g_0),\phi\circ\al_1,\phi(g_1),\ldots, \phi\circ\al_n, \phi(g_n)]$$ on arrows.

Let $i_{\cK}: \cK\to \cK_*$ be the identity map on objects and $i_{\cK}(g)=[g]$ on arrows. We have the following commutative diagram of morphisms of Lie groupoids:
$$         
\xymatrix{
\cK \ar^{i_{\cK}}[d] \ar^{\phi}[r] & \G \ar^{i_{\G}}[d] \\ 
\cK_* \ar^{\phi_*}[r] & \G_* } 
$$
 Haefliger's fundamental group of a groupoid \cite{Bri,H1,H2,H3} at $x_0$ , $\pi_1(\G,x_0)$,  coincides with the isotropy group ${(\G_*)}_{x_0}$ at $x_0$ of the fundamental groupoid $\G_*$.

\begin{prop}\cite{MM1}\label{ee}
\begin{enumerate}
\item If $\e: \cK\to \G$ is an essential equivalence, then $\e_*: \cK_*\to \G_*$ is an essential equivalence as well.
\item If $\cK\sim_M\G$ then $\cK_*\sim_M\G_*$ and the fundamental groups are isomorphic.
\item The fundamental groupoid $\G_{**}$ of $\G_*$ is isomorphic to $\G_*$.
\end{enumerate}
\end{prop}

\subsection{The bicategory $\HH$}\label{we}
Consider the category of Lie groupoids and functors. We are now ready to introduce a notion of  strict 1-homotopy  between functors.
\begin{definition} The morphisms $\phi:\cK\to \G$ and $\psi:\cK\to \G$ are {\it 1-homotopic} 
if their induced morphisms $\phi_*$ and $\psi_*$ between the fundamental groupoids are equivalent by a natural transformation. We write $\phi\simeq_{H}\psi$ .
\end{definition}
Since a natural transformation from  $\phi_*$ to $\psi_*$ associates to each object $x$ in $(G_*)_0=G_0$ an arrow $g_x=[g_0,\al_1,g_1,\ldots, \al_n, g_n]$ in $(G_*)_1$ from  $\phi(x)$ to $\psi(x)$, 
this  notion of homotopy corresponds to the intuitive idea of continuously deforming $\phi$ into $\psi$ by morphisms from $\cK$ to $\G$ along $\G$-paths.

We define a 2-morphism $H_a \colon\phi\ri\psi$ as a natural transformation $a\colon \phi_*\to\psi_*$:
\[
 \xy
 (0,-10)*{\bullet}="1"+(-3,0)*{\scriptstyle \cK_*};
 (14,-10)*{\bullet}="2"+(+3,0)*{\scriptstyle \G_*};
  (0,14)*{\bullet}="3"+(-3,0)*{\scriptstyle \cK};
   (14,14)*{\bullet}="4"+(+3,0)*{\scriptstyle \G};
  "1";"2" **\crv{(3,-1) & (11,-1)};
  "1";"2" **\crv{(3,-19) & (11,-19)};
   "3";"4" **\crv{(3,23) & (11,23)};
  "3";"4" **\crv{(3,5) & (11,5)};
(7,-3.25)*{\scriptstyle >}+(0,3)*{\scriptstyle \phi_*};
    (7,-16.75)*{\scriptstyle >}+(0,-3)*{\scriptstyle \psi_*};
    (7,20.75)*{\scriptstyle >}+(0,3)*{\scriptstyle \phi};
    (7,7.25)*{\scriptstyle >}+(0,-3)*{\scriptstyle \psi};
  {\ar@{=>}_{\scriptstyle a}(7,-6)*{};(7,-14)*{}};
   {\ar@{=>}_{\scriptstyle H_a}(7,18)*{};(7,10)*{}};
 "3"; "1" **\dir{-}? *\dir{>};
    "4"; "2" **\dir{-} ?*\dir{>};
\endxy
 \]
Lie groupoids, functors and 1-homotopies $H:\phi\ri\psi$ form a bicategory $\HH$.
All the 2-morphisms in $\HH$ are isomorphisms.

Horizontal and vertical compositions of 2-morphisms are given by the horizontal and vertical compositions of natural transformations, $a_*b_*$ and $a_*\cdot b_*$ respectively.

This notion of homotopy generalizes the concepts of natural transformation and ordinary homotopy, we have the following
\begin{prop} \label{nt} Let $\phi:\cK\to \G$ and $\psi:\cK\to \G$ be morphisms of  Lie groupoids.
\begin{enumerate}
\item If $\phi\sim_T\psi$ where $T$ is a natural transformation, then there is a 2-morphism $H:\phi\ri\psi$ in $\HH$.
\item If $\phi\simeq_F\psi$ where $F$ is an ordinary homotopy, then there is a 2-morphism $H:\phi\ri\psi$ in $\HH$.
\end{enumerate}
\end{prop}
\begin{proof}  
We will construct in each case a natural transformation $a:K_0\to G_{1*}$ satisfying $\psi([g_0,\al_1,g_1,\ldots, \al_n, g_n])a(x)=a(y)\phi([g_0,\al_1,g_1,\ldots, \al_n, g_n]$.
\begin{enumerate}
\item If $T:K_0\to G_1$ is a natural transformation with $T(x) : \phi(x)\to \psi(x)$ an arrow in $G_1$, define a natural transformation $a:K_0\to G_{1*}$ by 
$a(x)=[T(x)]$. We have that $a(x)$ is an arrow in $G_{1*}$ from $s(T(x))=\phi(x)$ to $t(T(x))=\psi(x)$ verifying the required equality.
\item An ordinary homotopy $F:\cK\times I\to \G$ with $F_0=\phi$ and $F_1=\psi$ determines for each $x\in K_0$ a path $F_x:I\to G_0$ from $\phi(x)$ to $\psi(x)$. Define 
$a:K_0\to G_{1*}$ by $a(x)=[1_{\phi(x)}, F_x, 1_{\psi(x)}]$
\[
 \xy
 (13,0)*{\bullet}="1";
 (27,0)*{\bullet}="2";
 "1";"1" **\crv{(0,9)& (26,9)}?(.55)*\dir{>};
  "2";"2" **\crv{(14,9)& (40,9)}?(.55)*\dir{>};
       "1"; "2" **\dir{~} *\dir{>};
        (20,2.5)*{{\scriptstyle F_x}};
 \endxy 
 \]

\end{enumerate}
\end{proof}

A strict {\it 1-homotopy equivalence} is a morphism $\phi:\cK\to \G$ such that there exists another morphism $\varphi:\G\to \cK$ and 2-isomorphisms  $\phi\varphi\ri \id_{\G}$  and $\varphi\phi\ri \id_{\cK}$ in $\HH$. We will say that two groupoids $\cK$ and $\G$ have the same strict {\it 1-homotopy type} if they are equivalent in the bicategory $\HH$.
But this notion of homotopy {\it is not} invariant of Morita equivalence. We need to add more morphisms and 2-morphisms to the bicategory $\HH$ and define our notion of Morita homotopy in an extended bicategory.

We can  characterize the strict 1-homotopy equivalences as  the morphisms that induce an equivalence between the fundamental groupoids. Recall that $\GG$ is the 2-category of Lie groupoids, functors and natural transformations.

\begin{prop}\label{se}
If $\phi:\cK\to \G$ is a 1-homotopy equivalence in $\HH$, then $\phi_*:\cK_*\to \G_*$ is an equivalence in $\GG$.
\end{prop}
\begin{proof}
We have the following diagram:
\[
 \xy
  (-15,0)*{\bullet}="1";
  (0,9)*{\bullet}="2";
  (15,0)*{\bullet}="3";
 "1";"2" **\crv{(-12,9)}?(.55)*\dir{>};
   "2";"3" **\crv{(12,9)}?(.55)*\dir{>};
   "1"; "3" **\dir{-} *\dir{>}; 
    (-10,6.5)*{}+(0,3)*{\scriptstyle \phi};
    (10,6.5)*{}+(0,3)*{\scriptstyle \varphi};
          (-17.5,0)*{\scriptstyle \cK};
  (17.5,0)*{\scriptstyle \cK};
  {\ar@{=>}_{\scriptstyle H}(0,7)*{};(0,1)*{}} ;
    (0,-3)*{\scriptstyle \id_{\cK}};
  \endxy .
 \]
then $(\varphi\phi)_*\sim_a (\id_{\cK})_*$ where $a: K_0\to (K_*)_1$ is a natural transformation. Since $(\varphi\phi)_*=\varphi_*\phi_*$ and $(\id_{\cK})_*=\id_{\cK_*}$, we have that 
$\varphi_*\phi_*\sim_a \id_{\cK_*}$. In the same way  $\phi_*\varphi_*\sim_b \id_{\G_*}$ and $\phi_*$ is an equivalence.
\end{proof}

Consider a 2-functor $\pi: \GG \to \GG$ between 2-categories given by $\pi(\G)=\G_*$, $\pi(\phi)=\phi_*$ and $\pi(T)=T_*$, where $T_*:\phi_*\ri\psi_*$ is a natural transformation defined in the following way. For each $x\in (K_*)_0=K_0$, we define $T_*(x):\phi(x)\to\psi(x)$ as the arrow in $\G_*$ given by $T_*(x)=[T(x)]$. This arrow satisfies the equality 
$\psi([g_0,\al_1,g_1,\ldots, \al_n, g_n])T_*(x)=T_*(y)\phi([g_0,\al_1,g_1,\ldots, \al_n, g_n])$. Then two morphisms $\phi$ and $\psi$ are 1-homotopic if their images by $\pi$ are equivalent. All the equivalences in $\GG$ are 1-homotopy equivalences in $\HH$, but there are 1-homotopy equivalences that do not come from  equivalences in $\GG$.

Now we introduce the {\it essential 1-homotopy equivalences} as the morphisms that induce an essential equivalence between the fundamental groupoids.
\begin{definition}
A morphism $\phi:\cK\to \G$ is an essential 1-homotopy equivalence if $\phi_*:\cK_*\to \G_*$ is an essential equivalence.
\end{definition}
In this case, $\phi_*$ defines an isomorphism between fundamental groups.

Let $E$ be the set of essential equivalences in $\GG$ and $W$ the set of essential 1-homotopy equivalences in $\HH$. We will show the following implications:
$$
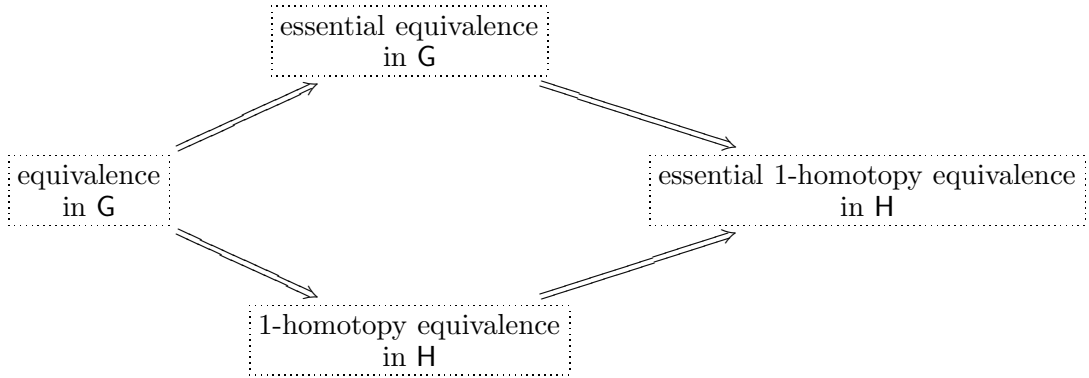
\xymatrix{
 & \begin{xy} 
*+{\txt{essential equivalence \\%
in $\GG$}}*\frm{.} 
\end{xy}\ar@2{->}[rd] & \cr
\begin{xy} 
*+{\txt{equivalence \\%
in $\GG$}}*\frm{.} 
\end{xy}\ar@2{->}[ur]\ar@2{->}[dr] & & \begin{xy} 
*+{\txt{essential 1-homotopy equivalence \\%
in $\HH$}}*\frm{.} 
\end{xy} \cr
 & \begin{xy} 
*+{\txt{1-homotopy equivalence\\%
in $\HH$}}*\frm{.} 
\end{xy}\ar@2{->}[ru] & 
}
$$

We have that every equivalence in $\GG$ is an essential equivalence in $\GG$ \cite{MM}. We will see that it is also a 1-homotopy equivalence 
in $\HH$.
\begin{prop} If $\phi:\cK\to \G$ is an equivalence in $\GG$, then $\phi$ is a 1-homotopy equivalence 
in $\HH$.
\end{prop}
\begin{proof} Let $\psi:\G\to \cK$ be the inverse up to equivalence in $\GG$. Then $\phi\psi\sim_T\id_{\G}$ and $\psi\phi\sim_{T'}\id_{\cK}$. By Proposition \ref{nt} (1) there are 2-morphisms $H:\phi\psi\ri\id_{\G}$  and $H':\psi\phi\ri\id_{\cK}$ in $\HH$.

\end{proof}

\begin{prop} \label{BF1}
\begin{enumerate}
\item Every essential equivalence in $\GG$  is an essential 1-homotopy equivalence in $\HH$.
\item  In $\HH$, all  1-homotopy equivalences are essential 1-homotopy equivalences.
\end{enumerate}
\end{prop}
\begin{proof} 
If $\e:\cK\to \G$ is an essential equivalence, then $\e$ is an essential 1-homotopy equivalence since it induces an essential equivalence between fundamental groupoids by Proposition \ref{ee}(1).

If $\phi:\cK\to \G$ is a 1-homotopy equivalence in $\HH$, then $\phi_*:\cK_*\to \G_*$ is an equivalence in $\GG$. Then  $\phi_*$ is also an essential equivalence in $\GG$ and $\phi$ is an essential 1-homotopy equivalence in $\HH$.
\end{proof}

\begin{lemma}\label{tran} For all $\e:\cJ\to\G$ and $\phi:\cK\to\G$ with $\e$ being an essential 1-homotopy equivalence, there exists an object $\cP$ and morphisms $\delta:\cP\to\cK$ and $\psi:\cP\to\cJ$ with $\delta$ essential 1-homotopy equivalence such that the following square commutes up to a 2-isomorphism:
$$         
\xymatrix{
\cP \ar^{\delta}[d] \ar^{\psi}[r] & \cJ\ar^{\e}[d] \\ 
\cK \ar^{\phi}[r] & \G } 
$$

\end{lemma}
\begin{proof}
We start by defining the {\it weak homotopy pullback} $\cP$ of the morphisms 
$$         
\xymatrix{
& \cJ\ar^{\e}[d] \\ 
\cK \ar^{\phi}[r] & \G} 
$$
as follows. Objects are triples $(x, [g_0,\al_1,\ldots, \al_n, g_n], y)$ where $x\in J_0$, $y\in K_0$ and $[g_0,\al_1,\ldots, \al_n, g_n]$ is a $\G$-path from $\e(x)$ to $\phi(y)$. Arrows in $\cP$ from 
$(x, [g_0,\al_1,\ldots, \al_n, g_n], y)$ to $(x', [g'_0,\al'_1,\ldots, \al'_n, g'_n], y')$ are pairs $(j,k)$ of arrows $j\in J_1$ and $k\in K_1$ such that
$$[g'_0,\al'_1,\ldots, \al'_n, g'_n][\e(j)]=[\phi(k)] [g_0,\al_1,\ldots, \al_n, g_n]$$
We observe that $\cP$ is the ordinary weak pullback $\cJ\times_{\G_*}\cK$ of the morphisms $\phi_*i_{\cK}$ and $\e_*i_{\cJ}$. Since $\e_*$ is an essential equivalence and $i_{\cJ}$ is the identity on objects we can assume that $\e_*i_{\cJ}$ is a submersion on objects and $\cP$ is a Lie groupoid. The square 
$$         
\xymatrix{
\cP\ar^{p_1}[r] \ar^{p_3}[d] & \cJ\ar^{\e}[d] \\ 
\cK \ar^{\phi}[r] & \G} 
$$
does not necessarily commute but it does when taken the induced morphisms in the fundamental groupoids. Consider the weak pullback $\cJ_*\times_{\G_*}\cK_*$ of groupoids 
$$         
\xymatrix{
\cJ_*\times_{\G_*}\cK_*\ar^{\;\;\pi_1}[r] \ar^{\pi_3}[d] & \cJ_*\ar^{\e_*}[d] \\ 
\cK_* \ar^{\phi_*}[r] & \G_*} 
$$
where $\e_*\pi_1\sim \phi_*\pi_3$ and $\pi_3$ is an essential equivalence. By definition of weak pullback and the explicit description of arrows in $\cP$ and $\cP^*$, we have that $\cJ_*\times_{\G_*}\cK_*=\cP_*$, $\pi_1 ={p_1}_*$ and $\pi_3 ={p_3}_*$, then $\cP_*$ is the weak pullback of $\e_*$ and $\phi_*$. Since the weak pullback square commutes up to natural transformation, we have that there is a 2-morphism $H:\e p_1\ri\phi p_3$ with $p_3\in W$.

\end{proof}

\begin{definition} Two Lie groupoids $\cK$ and $\G$ are  {\it Morita 1-homotopy equivalent} if there exist essential 1-homotopy equivalences:
\[\cK\overset{\w}{\gets}\cL\overset{\theta}{\to}\G\] for a third  Lie groupoid $\cL$.
\end{definition}
This defines an equivalence relation that we denote  $\simeq_M$. The transitivity property follows from Lemma \ref{tran}.
We will show next that the set of essential 1-homotopy equivalences $W$ admits a bicalculus of right fractions in $\HH$.

\subsection{Bicategories of fractions} Given a bicategory $\BB$ and a subset $S\subset B_1$ of morphisms satisfying certain conditions, there exists a bicategory $\BB[S^{-1}]$ having the same objects as $\BB$ but inverse morphisms of morphisms in $S$ have been added as well as more 2-morphisms. This new bicategory is called a bicategory of fractions of $\BB$ with respect to $S$ and was constructed by Pronk in \cite{P}. The bicategory of fractions $\BB[S^{-1}]$ is characterized by the universal property that any homomorphism $F:\BB\to\DD$ sending elements of $S$ into equivalences factors as $F\cong\tilde F\circ U$
$$         
\xymatrix{
\BB \ar^{F}[d] \ar^{U\;\;\;\;}[r] & \BB[S^{-1}]\ar^{\tilde F}[ld] \\ 
\DD & } 
$$
where $\tilde F$ is unique with the property up to isomorphism.
The following conditions are needed on $S$ to admit a bicalculus of right fractions  \cite{P}:
\begin{enumerate}
\item[BF1] All equivalences are in $S$.
\item[BF2] If $\phi$ and $\psi$ are in $S$, then $\phi\psi\in S$.
\item[BF3] For all $\e:\cJ\to\G$ and $\phi:\cK\to\G$ with $\e\in S$ there exists an object $\cP$ and morphisms $\delta:\cP\to\cK$ and $\psi:\cP\to\cJ$ with $\delta\in S$ such that the following square commutes up to a 2-isomorphism:
$$         
\xymatrix{
\cP \ar^{\delta}[d] \ar^{\psi}[r] & \cJ\ar^{\e}[d] \\ 
\cK \ar^{\phi}[r] & \G } 
$$
\item[BF4] If $H:\eta\phi\ri\eta\varphi$ is a 2-morphism with $\eta\in S$ then there exists a morphism $\e\in S$ and a 2-morphism $G:\phi\e\ri\varphi\e$ such that $H\cdot\e=\eta\cdot G$. When $H$ is an isomorphism, we require $G$ to be an isomorphism too. When $\e'$ and $G'$ are another such pair, there exist 1-morphisms $\delta$, $\delta'$ such that $\e\delta$ and 
$\e'\delta'$ are in $S$ and a 2-isomorphism $J\colon \e\delta\ri\e'\delta'$  such that the following diagram commutes:
$$         
\xymatrix{
\phi\e\delta \ar^{\phi J}[d] \ar^{G\delta}[r] & \psi\e\delta\ar^{\psi J}[d] \\ 
\phi\e'\delta' \ar^{G'\delta'}[r] & \psi\e'\delta'.}
$$

\item[BF5] If there is a 2-isomorphism $\e\ri\delta$ with $\delta\in S$ then $\e\in S$.
\end{enumerate}
We will prove now that the set of essential 1-homotopy equivalences $W$ allows a bicalculus of right fractions.

Conditions BF1 and BF2 follow from Proposition \ref{BF1} and from the fact that the set of essential equivalences is closed under composition. BF3 is Lemma \ref{tran}.
If $\phi, \psi: \cK\to \G$ are morphisms, $\eta:\G\to \cJ$ is an essential 1-homotopy equivalence and $H:\eta\phi\ri\eta\varphi$ is a 2-morphism, we have that $\eta_*$ is an essential equivalence. Then there exist a natural transformation $b$ \cite{M1} with $\phi_*\sim_b\psi_*$ and $\eta_* b=a$, where $a$ is the natural transformation between $\eta_*\phi_*$ and $\eta_*\varphi_*$. The natural transformation $b$ defines a 2-morphisms $G:\phi\ri\psi$ such that $H=\eta G$ and condition BF4 follows.

Finally, if $\e\ri\delta$ is a 2-isomorphism, then there is a natural transformation between $\e_*$ and $\delta_*$. Since $\delta$ is an essential 1-homotopy equivalence, $\delta_*$ is an essential equivalence and this implies that $\e_*$  is an essential equivalence as well. Then, $\e$ is an essential 1-homotopy equivalence.

Therefore, there exists a bicategory of fractions $\HH(W^{-1})$ inverting the essential 1-homotopy equivalences.

\subsection{The bicategory $\HH(W^{-1})$} The objects of $\HH(W^{-1})$ are Lie groupoids. The morphisms from $\cK$ to $\G$ are formed by pairs $(\w,\phi)$
\[\cK\overset{\w}{\gets}\cJ\overset{\phi}{\to}\G\]
such that $\w$ is an essential 1-homotopy equivalence. 

The composition of morphisms $(\G\overset{\w'}{\gets}\cJ'\overset{\phi'}{\to}\cL)\circ (\cK\overset{\w}{\gets}\cJ\overset{\phi}{\to}\G)$  is given by a morphism
$$\cK\xleftarrow{\w p_3}\cP\xrightarrow{\phi' p_1}\cL$$
where  $\cP$ is  the weak homotopy pullback of $\w'$ and $\phi$.

A 2-morphism from $(\w,\phi)$ to $(\w',\phi')$ is given by a class of  diagrams:
$$
\xymatrix{ & 
{\cJ}\ar[dr]^{\phi}="0" \ar[dl]_{\w}="2"&\\
{\cK}&{\cL} \ar[u]_{u} \ar[d]^{v}&{\G}\\
&{\cJ'}\ar[ul]^{\w'}="3" \ar[ur]_{\phi'}="1"&
\ar@{}"0";"1"|(.4){\,}="7" 
\ar@{}"0";"1"|(.6){\,}="8" 
\ar@{=>}"7" ;"8"_{H'} 
\ar@{}"2";"3"|(.4){\,}="5" 
\ar@{}"2";"3"|(.6){\,}="6" 
\ar@{=>}"5" ;"6"^{H} 
}
$$
where $\cL$ is a Lie groupoid, $u$ and $v$ are essential 1-homotopy equivalences and $H:\w u\ri \w' v$ and $H':\phi u\ri \phi'v$ are 2-isomorphisms in $\HH$. The horizontal and vertical composition of diagrams are defined as in \cite{P}.

The notion of 1-homotopy we propose corresponds to 2-morphisms in the bicategory $\HH(W^{-1})$. That is, we will say that two morphisms  are  1-homotopic if there is a 2-morphism between them:

\begin{definition}
Two morphisms $(\w,\phi)$ and $(\w',\phi')$ are {\it 1-homotopic} if there is there exists a diagram 
$$
\xymatrix{ & 
{\cJ}\ar[dr]^{\phi}="0" \ar[dl]_{\w}="2"&\\
{\cK}&{\cL} \ar[u]_{u} \ar[d]^{v}&{\G}\\
&{\cJ'}\ar[ul]^{\w'}="3" \ar[ur]_{\phi'}="1"&
\ar@{}"0";"1"|(.4){\,}="7" 
\ar@{}"0";"1"|(.6){\,}="8" 
\ar@{=>}"7" ;"8"_{H'} 
\ar@{}"2";"3"|(.4){\,}="5" 
\ar@{}"2";"3"|(.6){\,}="6" 
\ar@{=>}"5" ;"6"^{H} 
}
$$
as above.
 \end{definition}
In this case, we write $(\w,\phi)\simeq(\w',\phi')$ and we say that there is a {\it 1-homotopy} between $(\w,\phi)$ and $(\w',\phi')$.

In particular, when $\w$ and $\w'$ are essential equivalences, we have a notion of 1-homotopy for generalized maps and when they are identities, we have a notion of 1-homotopy for strict maps.

This notion of 1-homotopy relates to the Morita homotopy notion defined in \cite{Colman}.
\begin{remark}\cite{Colman}
A Morita homotopy equivalence $\cK\simeq_M\G$ induces a Morita equivalence between the homotopy groupoids $\cK^*\sim_M\G^*$. If $\cK$ and $\G$ are Morita homotopy equivalent, then they are $1$-homotopy equivalent.
\end{remark}

Two objects $\cK$ and $\G$  are equivalent in $\HH(W^{-1})$  if there are morphisms $(\w,\phi)$ from  $\cK$ to $\G$ and $(\theta,\psi)$ from $\G$ to $\cK$ such that $(\w,\phi)\circ(\theta,\psi)$ is 1-homotopic to the identity $(\id_{\G},\id_{\G})$ and $(\theta,\psi)\circ(\w,\phi)\simeq(\id_{\cK},\id_{\cK})$.

\begin{prop} A morphism $(\w,\phi)$ is invertible up to a 2-isomorphism  in $ \HH(W^{-1})$ if and only if $\phi$ is an essential 1-homotopy equivalence. In this case, the inverse of $(\w,\phi)$ is the morphism $(\phi,\w)$.
\end{prop}
 In other words, the definition of Morita 1-homotopy equivalence in subsection \ref{we} amounts to {\it equivalence} in the bicategory $\HH(W^{-1})$. So, we write $\cK\simeq_M\G$ for equivalence of objects in $\HH(W^{-1})$. The {\it 1-homotopy type}  of $\G$ is the class of $\G$ under the equivalence relation $\simeq_M$.

We show now that the 1-homotopy type is invariant under Morita equivalence.
\begin{prop} If $\cK\sim_M\G$,  then $\cK\simeq_M\G$.
\end{prop}
\begin{proof}
If $\cK$ and $\G$ are Morita equivalent, then there is a Lie groupoid $\cL$ and essential equivalences:
\[\cK\overset{\e}{\gets}\cL\overset{\delta}{\to}\G.\]
The maps $\e$ and $\delta$ are also essential 1-homotopy equivalences by Proposition \ref{BF1} (1). Then the morphisms $(\e,\delta)$ and $(\delta,\e)$ are inverse up to a 2-isomorphism in $\HH$. Then $\cK$ is equivalent to $\G$ in the bicategory $\HH(W^{-1})$.
\end{proof}

The bicategory $\HH(W^{-1})$ comes equipped with a universal homomorphism  \cite{P}  $U:\HH\to \HH(W^{-1})$  sending essential 1-homotopy equivalences to equivalences and with the following universal property:  any homomorphism of bicategories $F:\HH\to \BB$ sending essential 1-homotopy equivalences to equivalences factors in a unique way as $F=G\circ U$ where  $G:\HH(W^{-1})\to \BB$ is a homomorphism of bicategories.

We have the following square of bicategories which relates in particular the notions of equivalence and Morita equivalence in $\GG$ and the 1-homotopy equivalence and Morita 1-homotopy equivalence in $\HH$:
$$
\xymatrix{ 
\GG
\ar[r]^{U|_{\GG}} 
\ar[d]_{i_{\GG}} 
& \GG(E^{-1})\ar[d]^{i_{\GG(E^{-1})}}\\ 
\HH
\ar[r]^{U} 
& \HH(W^{-1})
}
$$
This diagram commutes on the nose by choosing a suitable $i_{\GG(E^{-1})}$.
\subsection{The classifying space}
We will show that the homotopy groups of a groupoid defined in terms of the classifying space are invariant of the Morita 1-homotopy type as defined above. First we recall the construction of the homotopy groups of a groupoid as the homotopy groups of its classifying space \cite{Se,H1}.

The classifying space $B\G$ of a groupoid $\G$ is the geometric realization of the simplicial manifold whose n-simplices are the contravariant functors $[n]\to \G$ where $[n]$ is the linearly ordered set $\{ 0,1,\ldots, n\}$. We can describe the simplicial manifold $G_*$ as a sequence of manifolds of composable strings of arrows:
$$G_n=\{(g_1,\dots, g_n) \; |\; g_i\in G_1, s(g_i)=t(g_{i+1}), i=1,\ldots, n-1\}$$
connected by face operators
$$d_i:G_n\to G_{n-1}$$ given by 
$d_i(g_1,\dots, g_n) =
\begin{cases} (g_2,\dots, g_n) & i=0
\\
(g_1,\dots,  g_{i}g_{i+1}, \dots,   g_{n}) & 0<i<n
\\
(g_1,\dots, g_{n-1}) &i=n
\end{cases}
$

when $n>1$ and $d_0(g)=s(g)$, $d_1(g)=t(g)$ when $n=1$.

The classifying space of $\G$ is the geometric realization of the simplicial manifold $G_*$:

$$B\G=|G_*|=\bigsqcup_n(G_n\times\Delta^n)/ (d_i(g),x)\sim (g, \delta_i(x))$$
where $\Delta^n$ is the standard topological $n$-simplex and $\delta_i:\Delta^{n-1}\to \Delta^n$ is the linear embedding of the $i$-th face.

\begin{definition} The fundamental group of a groupoid $\G$ is defined as the fundamental group of $B\G$: $$\pi_1(\G, x)=\pi_1(B\G, x)$$
where $x\in G_0$.
\end{definition}

\begin{remark}
This definition of the fundamental group of a groupoid $\G$ by Segal coincides with Haefliger's definition presented in section 3 \cite{Se1}.
\end{remark}

\begin{prop} \cite{Se,M1} A morphism $\phi:\cK\to \G$ induces a map $B\phi:B\cK\to B\G$. If $\phi$ is an essential equivalence, then $B\phi$ is a weak homotopy equivalence.
\end{prop}
If two groupoids are Morita equivalent, then they have isomorphic homotopy groups. We will show that this remains true for the much weaker invariant of Morita 1-homotopy type.

\begin{prop} If $\phi:\cK\to \G$ is an essential 1-homotopy equivalence, then $\pi_1(\cK, x)=\pi_1(\G, \phi(x))$
\end{prop}
\begin{proof} The morphism $\phi:\cK\to \G$  induces an essential equivalence $\phi_*:\cK_*\to \G_*$ on the fundamental groupoids. Then $B\phi_*:B\cK_*\to B\G_*$ is a weak homotopy equivalence and we have that $\pi_1(B\cK_*, x)=\pi_1(B\G_*, \phi(x))$. Since $\pi_1(\cK_*, x)=\pi_1(B\cK_*, x)$ by definition and $\pi_1(\cK_*)=\cK_*$ by Proposition \ref{ee} (3), we have that $\pi_1(B\cK_*, x)=\pi_1(\cK, x)$. Analogously, 
$\pi_1(B\G_*, \phi(x))=\pi_1(\G, \phi(x))$ and the claim follows.
\end{proof}

\section{Orbifolds as groupoids}
We recall now the description of orbifolds as groupoids due to Moerdijk and Pronk \cite{MP,M}. Orbifolds were first introduced by Satake \cite{S} as a generalization of a manifold defined in terms of local quotients. The groupoid approach provides a global language to reformulate the notion of orbifold.

We follow the exposition in \cite{A}. A groupoid $\G$ is {\it proper} if $(s,t):G_1\to G_0\times G_0$ is a proper map and it is a {\it foliation} groupoid is each isotropy group is discrete. 

\begin{definition}
An {\it orbifold} groupoid is a proper foliation groupoid.
\end{definition}

For instance the holonomy group of a foliation $\F$ is always a foliation groupoid but it is an orbifold groupoid if and only if $\F$ is a compact-Hausdorff foliation.

Given an orbifold groupoid $\G$, its orbit space $|\G|$ is a locally compact Hausdorff space. Given an arbitrary locally compact Hausdorff space $X$ we can equip it  with an orbifold structure as follows:

\begin{definition} An {\it orbifold structure} on a locally compact Hausdorff space $X$ is given by an orbifold groupoid $\G$ and a homeomorphism $f:|\G|\to X$.
\end{definition}
If $\e:\cK\to \G$ is an essential equivalence and $|\e|:|\cK|\to |\G|$ is the induced homeomorphism between orbit spaces, we say that the composition $f\circ|\e|:|\cK|\to X$ defines an {\it equivalent} orbifold structure. 

\begin{definition} An {\it orbifold}  $\X$ 
is a space $X$ equipped with an equivalence class of orbifold 
structures. A specific such structure, given by 
$\G$ and  $f : |\G | \to X $ is
a {\it presentation} of the orbifold 
$\X$.
\end{definition}

If two groupoids are Morita equivalent, then they define the same orbifold. Therefore any structure or invariant for orbifolds, if defined through groupoids, should be invariant under Morita equivalence. 

An orbifold map $\Y\to \X$ is defined as a class of generalized maps $(\e,\phi)$ from $\cK$ to $\G$ between presentations of the orbifolds such that the diagram commutes:

\[\xymatrix{
|\cK| \ar[r] \ar[d]& |\G| \ar[d]^{f} \\ 
Y \ar[r]& X}\]
A specific such generalized map $\cK\overset{\e}{\gets}\cJ\overset{\phi}{\to}\G$ 
is called a {\it presentation} of the orbifold map.

Therefore, our Morita invariant notion of 1-homotopy for generalized maps gives a notion of 1-homotopy for orbifolds:
\begin{definition}
Two orbifold maps $\xymatrix{\Y \ar@<.5ex>[r] \ar@<-.5ex>[r] & \X} $ are {\it 1-homotopic} if their presentations $(\e,\phi)$ and $(\nu,\psi)$ are 1-homotopic.
\end{definition}
In other words, if $\cK\overset{\e}{\gets}\cJ\overset{\phi}{\to}\G$ and $\cK\overset{\nu}{\gets}\cJ'\overset{\psi}{\to}\G$  are presentations of  the orbifold maps, they are 1-homotopic if there exists a Lie groupoid $\cL$ and essential 1-homotopy equivalences 
$$\cJ\overset{\w}{\gets}\cL\overset{\w'}{\to}\cJ'$$
such that $\e_*\w_*\sim\nu_*\w'_*$ and $\phi_*\w_*\sim\psi_*\w'_*$.

The equivalence class of a morphism in $\HH(W^{-1})$ defines the 1-homotopy class of the orbifold map. If $\cK$ and $\G$ are presentations of the orbifold $\Y$ and $\X$ respectively, a specific morphism $\cK\overset{\w}{\gets}\cJ\overset{\phi}{\to}\G$ with $\w$ being an essential 1-homotopy equivalence, is a presentation of the 1-homotopy class of $\Y\to \X$.

With the obvious notion of {\it orbifold 1-homotopy type}, we have that 
the orbifold fundamental group of $\X$ defined as the fundamental group of a presentation groupoid $\G$ \cite{H1,H2,H3,C}, $$\pi_1^{\mbox{\tiny\it  orb}}(\X, \bar x)=\pi_1(\G, x)$$ is an invariant of orbifold 1-homotopy type.

\begin{example}
Consider the orbifold $\X$ having as a presentation groupoid the holonomy groupoid associated to the Seifert fibration on the M\"{o}bius band, $\G=\Hol(M,\F_S)$ and the orbifold $\Y$ represented by the point groupoid $\bullet^{\Z_2}$. These groupoids are not Morita equivalent and they do not have the same type of (strict) homotopy.
Consider the morphism $$\bullet^{\Z_2}\overset{c}{\gets}\cJ\overset{i}{\to}\G$$ where 
$\cJ=\Hol_I(M,\F_S)$ is the reduced holonomy groupoid to a transversal interval $I$,  and $c$ is the constant map on objects $J_0=I$ and the constant map on each connected component of the manifold of arrows $J_1=I\sqcup I$. This morphism $(i,c)$  is invertible since both the constant map $c$ and the inclusion map $i$ are essential 1-homotopy equivalences. The inverse morphism is given by
$$\G\overset{i}{\gets}\cJ\overset{c}{\to}\bullet^{\Z_2}$$
and the orbifolds $\X$  and $\Y$ have the same 1-homotopy type.
\end{example}

\end{document}